\documentclass[12pt, a4, ams]{article}
\usepackage{amsmath, amssymb}
\begin{document}
\newtheorem{Definition}{Definition}[section]
\newtheorem{Theorem}{Theorem}[section]
\newtheorem{Lemma}[Theorem]{Lemma}
\newtheorem{Corollary}[Theorem]{Corollary}

\title{Variance of the number of Comparisons of Randomised Quicksort}
\author{Vasileios Iliopoulos and David Penman\\
Department of Mathematical Sciences\\ University of Essex, CO4 3SQ, U.K\\
\texttt{viliop@essex.ac.uk, dbpenman@essex.ac.uk}}
\date{May 2010}
\maketitle

\tableofcontents

\pagebreak

\section{Introduction}

In this paper, we present a calculation of the variance of the number of comparisons required by the Quicksort algorithm for sorting a set, when the pivot is chosen uniformly and at random from the $n$ objects $\{x_{1},\ldots x_{n}\}$ (which have a total order on them, but not one initially known to us) to be sorted. Remember that, given a pivot $x_{i}$, the Quicksort proceeds by carrying out pairwise comparisons (which we assume can be done) of all the other objects with $x_{i}$, and using this to split the original set into two subsets, all those elements above the pivot and all those below it. We then iterate this process, choosing pivots in each smaller set uniformly at random and using comparisons with the pivot to split each set into two others. Eventually we will have all the elements in order and the algorithm terminates. The object of interest is the number $C_{n}$ of comparisons required to get the $n$ elements in order. If pivots in each set are chosen from all elements in the set uniformly at random, $C_{n}$ is clearly a random variable.
It is well-known that the mean $M_{n}$ of $C_{n}$ is equal to $2(n+1)H_{n}-4n$, where $H_{n}=1+1/2+1/3+\ldots 1/n$ is the $n$th harmonic number. (Note that $H_{0}=0$).  For proofs of this fact, see \cite{prob}, \cite{fill}. We also define the $n$th harmonic number of order $k$ to be equal to $H^{(k)}_{n}=1+1/2^{k}+1/3^{k}+\ldots 1/n^{k}$.

In this paper,  we obtain the variance of $C_{n}$. The formula for this is stated without proof in Knuth \cite{knut}, who in his Exercise 6.2.2-8 states the formula
$$Var(C_{n})= 7n^{2}-4(n+1)^{2}H_{n}^{(2)}-2(n+1)H_{n}+13n.$$
Similarly, the papers \cite{nip} and \cite{sk} provide sketches of how to
prove this fact. Also, in \cite{reg} the
asymptotic variance of the random variable
$$\frac{C_{n}-2(n+1)(H_{n+1}-1)}{n+1}$$ is obtained using results about
moments of \lq the depth of insertion\rq\, in a tree and some martingale
arguments. However we are not aware of any source where all details of the
argument are written out explicitly with as few prerequisites as possible.
Thus we felt it would be desirable to provide such an account, though we
freely acknowledge that not all the details of the computation are
particularly interesting. No originality is claimed for the result. 

The basic strategy of the argument is to use a sequence of reductions of the problem. We first use generating functions to show that it is sufficient to prove that a certain sequence $B_{n}$, defined to the next section is equal to
$$B_{n}
=2(n+1)^{2}(H^{2}_{n}-H^{(2)}_{n})-H_{n}(n+1)(8n+2)+\frac{23n^{2}+17n}{2}.$$
The proof of this in turn relies on various identities involving harmonic numbers and much manipulative algebra  - readers may prefer to use MAPLE at some stages (as we did ourselves to initially find the relationships, though we do include proofs for completeness).

\section{Proof}

The theorem that we will prove at this paper is
\begin{Theorem}
If $C_{n}$ is the number of comparisons used by Quicksort
with a pivot chosen uniformly at random, then
$$\mbox{Var}(C_{n})=7n^{2}-4(n+1)^{2}H_{n}^{(2)}-2(n+1)H_{n}+13n.$$
\end{Theorem}
We start with a recurrence for the generating function of $C_{n}$, namely
$f_{n}(z)=\sum_{k=0}^{n(n-1)/2}P(C_{n}=k)z^{k}$.  We will use this to
reduce the proof of the theorem to proving a certain recurrence formula
for a quantity to be called $B_{n}$ (defined below).
\begin{Theorem}
In Random Quicksort of $n$ objects, the generating functions $f_{i}$ satisfy
$$f_{n}(z)=\frac{z^{n-1}}{n}\sum_{j=1}^{n}f_{j-1}(z)f_{n-j}(z).$$
\end{Theorem}
{\bf Proof.} Using the following equation
$$C_{n}=C_{U_{n}-1}+{C_{n-U_{n}}} +n-1$$
we have that
\begin{align*}
P(C_{n}=k)&=\frac{1}{n}\sum_{m=1}^{n}P(C_{n}=k\vert U_{n}=m)\\
&=\frac{1}{n}\sum_{m=1}^{n}\sum_{j=1}^{k-(n-1)}P(C_{m-1}=j)P(C_{n-m}=k-(n-1)-j).
\end{align*}
(We are using here the fact that $C_{m-1}$ and $C_{n-m}$ are independent).
Thus
\begin{eqnarray*}
P(C_{n}=k)z^{k}=\frac{1}{n}\sum_{m=1}^{n}\sum_{j=1}^{k-(n-1)}
P(C_{m-1}=j)z^{j}P(C_{n-m}=k-(n-1)-j)z^{k-(n-1)-j}z^{n-1}.\\
\end{eqnarray*}
Multiplying by $z^{k}$ and summing over $k$, so as to get the generating
function $f_{n}$ of $C_{n}$ on the left, we obtain
\begin{align*}
f_{n}(z)&=\frac{1}{n}\sum_{k=1}^{n-1+j}\sum_{m=1}^{n}\sum_{j=1}^{k-(n-1)}P(C_{m-1}=j)z^{j}P(C_{n-m}=k-(n-1)-j)z^{k-(n-1)-j}z^{n-1}\\
&=\frac{z^{n-1}}{n}\sum_{m=1}^{n}\sum_{j=1}^{k-(n-1)}P(C_{m-1}=j)z^{j}
\sum_{k=1}^{n-1+j}P(C_{n-m}=k-(n-1)-j)z^{k-(n-1)-j}.
\end{align*}
Thus
\begin{equation*}
f_{n}(z)=\frac{z^{n-1}}{n}\sum_{m=1}^{n}f_{m-1}(z)f_{n-m}(z)\,\,(*)
\end{equation*}
as required. $\bullet$

This of course leads to a recursion for the variance, using the
well-known link between variance of a random variable $X$ and its
generating function $f_{X}(z)$:
\begin{eqnarray*}
Var(X)=f_{X}^{\prime \prime}(1)+f_{X}^{\prime}(1)-(f_{X}^{\prime}(1))^{2}.
\end{eqnarray*}
We use this formula together with equation (*) above. For the first
derivative, 
\begin{align*}
f_{n}^{\prime}(z)=\frac{(n-1)z^{n-2}}{n}\sum_{j=1}^{n}f_{j-1}(z)f_{n-j}(z)
+\frac{z^{n-1}}{n}\sum_{j=1}^{n}f^{\prime}_{j-1}(z)
f_{n-j}(z)+\frac{z^{n-1}}{n}\sum_{j=1}^{n}f_{j-1}(z)f^{\prime}_{n-j}(z).
\end{align*}
From standard properties of generating functions, $E(C_{n})=f^{\prime}_{n}(1)$. Differentiating again we obtain
\begin{align*}
f''_{n}(z)&=\frac{(n-1)(n-2)z^{n-3}}{n}\sum_{j=1}^{n}f_{j-1}(z)f_{n-j}(z)+\frac{(n-1)z^{n-2}}{n}\sum_{j=1}^{n}f^{\prime}_{j-1}(z)f_{n-j}(z)\\
&+\frac{(n-1)z^{n-2}}{n}\sum_{j=1}^{n}f_{j-1}(z)f^{\prime}_{n-j}(z)+\frac{(n-1)z^{n-2}}{n}\sum_{j=1}^{n}f^{\prime}_{j-1}(z)f_{n-j}(z)\\
&+\frac{z^{n-1}}{n}\sum_{j=1}^{n}f''_{j-1}(z)f_{n-j}(z)+\frac{z^{n-1}}{n}\sum_{j=1}^{n}f^{\prime}_{j-1}(z)f^{\prime}_{n-j}(z)+\frac{(n-1)z^{n-2}}{n}\sum_{j=1}^{n}f_{j-1}(z)f^{\prime}_{n-j}(z)\\
&+\frac{z^{n-1}}{n}\sum_{j=1}^{n}f^{\prime}_{j-1}f^{\prime}_{n-j}(z)
+\frac{z^{n-1}}{n}\sum_{j=1}^{n}f_{j-1}(z)f''_{n-j}(z).
\end{align*}
Setting $z=1$, we have (see \cite{nip})
\begin{align*}
f''_{n}(1)&=(n-1)(n-2)+\frac{2}{n}(n-1)\sum_{j=1}^{n}M_{j-1}+\frac{2}{n}(n-1)\sum_{j=1}^{n}M_{n-j}\\
&+\frac{1}{n}\sum_{j=1}^{n}(f''_{j-1}(1)+f''_{n-j}(1))+\frac{2}{n}
\sum_{j=1}^{n}M_{j-1}M_{n-j}.
\end{align*}
where $M_{j-1}, M_{n-j}$ are $f^{\prime}_{j-1}(1)$, $f^{\prime}_{n-j}(1)$,
\textit{i.e.}
the mean number of comparisons to sort a set of $(j-1)$ \& $(n-j)$ elements
respectively. Setting $B_{n}=f_{n}^{"}(1)/2$, we obtain
\begin{align*}
2B_{n}&=(n-1)(n-2)+\frac{2(n-1)}{n}\sum_{j=1}^{n}M_{j-1}+\frac{2(n-1)}{n}\sum_{j=1}^{n}M_{n-j}+\frac{1}{n}\sum_{j=1}^{n}(2B_{j-1}+2B_{n-j})\\
&+\frac{2}{n}\sum_{j=1}^{n}M_{j-1}M_{n-j}.
\end{align*}
But now, noting that $\sum_{j=1}^{n}M_{j-1}=\sum_{j=1}^{n}M_{n-j}$, as
both sums are $M_{1}+\ldots+M_{n-1}$ (using the fact that $M_{0}=0$), and
similarly that $\sum_{j=1}^{n}B_{j-1}=\sum_{j=1}^{n}B_{n-j}$, we see this is
\begin{eqnarray}
B_{n}={n-1\choose 2}+\frac{2(n-1)}{n}\sum_{j=1}^{n}M_{j-1}+\frac{2}{n}\sum_{j=1}^{n}B_{j-1}
+\frac{1}{n}\sum_{j=1}^{n}M_{j-1}M_{n-j}.
\end{eqnarray}
What this argument has shown for us is the following - compare \cite{sk}
where it is also shown that this recurrence has to be solved, though
no details of how to solve it are given.
\begin{Theorem}
In order to prove Theorem 2.1, it is sufficient to show that the recurrence
equation (1) for $B_{n}$ is satisfied by \cite{sk},
$$B_{n}=2(n+1)^{2}H_{n}^{2}-(8n+2)(n+1)H_{n}
+\frac{n(23n+17)}{2}-2(n+1)^{2}H_{n}^{(2)}.$$
\end{Theorem}

{\bf Proof.} If we get this formula, we then have
\begin{align*}
\mbox{Var}(C_{n})&=f_{n}^{"}(1)+f_{n}^{'}(1)-(f_{n}^{'}(1))^{2}=2B_{n}+2(n+1)H_{n}-4n-[2(n+1)H_{n}-4n]^{2}\\
&=4(n+1)^{2}H_{n}^{2}-2(8n+2)(n+1)H_{n}+2\frac{n(23n+17)}{2}-4(n+1)^{2}H_{n}^{(2)}\\
&+2(n+1)H_{n}-4n-[2(n+1)H_{n}-4n]^{2}\\
&=4(n+1)^{2}H_{n}^{2}-2(8n+2)(n+1)H_{n}+n(23n+17)-4(n+1)^{2}H_{n}^{(2)}\\
&+2(n+1)H_{n}-4n-4(n+1)^{2}H^{2}_{n}+16n(n+1)H_{n}-16n^{2}\\
&=7n^{2}+13n-4(n+1)^{2}H_{n}^{(2)}-2(n+1)H_{n}[(8n+2)-1-8n]\\
&=7n^{2}-4(n+1)^{2}H_{n}^{(2)}-2(n+1)H_{n}+13n
\end{align*}
as required. $\bullet$
\section{Solution of the recurrence for $B_{n}$}
We have to solve the $B_{n}$ recurrence. For the sum of the $M_{j-1}$,
the expected numbers of comparisons, we have
\begin{equation*}
\sum_{j=1}^{n}M_{j-1}=\sum_{j=1}^{n}[2jH_{j-1}-4(j-1)]\\
=2\sum_{j=1}^{n}jH_{j-1}-4\sum_{j=1}^{n}(j-1).
\end{equation*}
For the computation of the first sum, a Lemma follows.
\begin{Lemma}
For $n \in \mathbf N$
\begin{eqnarray*}
\sum_{j=1}^{n}jH_{j-1}=\frac{n(n+1)H_{n+1}}{2}-\frac{n(n+5)}{4}\mbox{~and~}
\sum_{j=1}^{n}M_{j-1}=n(n+1)H_{n+1}-\frac{5n^{2}+n}{2}
\end{eqnarray*}
\end{Lemma}
{\bf Proof.}
By induction on $n$, the case $n=1$ being trivial. Suppose that it holds for
all $n\leq k$. Then for $n=k+1$ we have
\begin{align*}
\sum_{j=1}^{k+1}jH_{j-1}&=\sum_{j=1}^{k}jH_{j-1}+(k+1)H_{k}=\frac{k(k+1)H_{k+1}}{2}-\frac{k(k+5)}{4}+(k+1)H_{k}\\
&=\frac{k(k+1)H_{k+1}}{2}-\frac{k(k+5)}{4}+(k+1)(H_{k+1}-\frac{1}{k+1})\\
&=\frac{(k+1)(k+2)H_{k+1}}{2}-\frac{k(k+5)}{4}-1
\end{align*}
\begin{align*}
&=\frac{(k+1)(k+2)H_{k+2}}{2}-\frac{k+1}{2}-\frac{k(k+5)}{4}-1\\
&=\frac{(k+1)(k+2)H_{k+2}}{2}-\frac{k^{2}+7k+6}{4}\\
&=\frac{(k+1)(k+2)H_{k+2}}{2}-\frac{(k+1)(k+6)}{4}\\
&=\frac{n(n+1)H_{n+1}}{2}-\frac{n(n+5)}{4}.
\end{align*}
giving the first claim. The second claim follows recalling that
$\sum_{j=1}^{n}(j-1)=n(n-1)/2$. $\bullet$

Now, we will compute the term
\begin{eqnarray}
\sum_{j=1}^{n}M_{j-1}M_{n-j}.
\end{eqnarray}
We shall use three Lemmas in the proof.
\begin{Lemma}
For $n \in \mathbf N$, it holds that
\begin{eqnarray*}
\sum_{j=1}^{n}M_{j-1}M_{n-j}=4\sum_{j=1}^{n}jH_{j-1}(n-j+1)H_{n-j}-\frac{8}{3}n(n^{2}-1)H_{n+1}+\frac{44n}{9}(n^{2}-1)
\end{eqnarray*}
\end{Lemma}

{\bf Proof.} To do this, we will again use the formula obtained previously for
$M_{j}$. We have
\begin{align*}
\sum_{j=1}^{n}M_{j-1}M_{n-j}&=\sum_{j=1}^{n}[(2jH_{j-1}-4j+4)(2(n-j+1)H_{n-j}-4n+4j)]\\
&=4\sum_{j=1}^{n}jH_{j-1}(n-j+1)H_{n-j}-8n\sum_{j=1}^{n}jH_{j-1}+8\sum_{j=1}^{n}j^{2}H_{j-1}\\
&-8\sum_{j=1}^{n}j(n-j+1)H_{n-j}+16n\sum_{j=1}^{n}j-16\sum_{j=1}^{n}j^{2}+8\sum_{j=1}^{n}(n-j+1)H_{n-j}\\
&-16n^{2}+16\sum_{j=1}^{n}j.
\end{align*}
We need to work out the value of $\sum_{j=1}^{n}j^{2}H_{j-1}$. Using
MAPLE initially, we found
\begin{eqnarray*}
\sum_{j=1}^{n}j^{2}H_{j-1}=\frac{6n(n+1)(2n+1)H_{n+1}-n(n+1)(4n+23)}{36};
\end{eqnarray*}
we will confirm this by induction.
\begin{Lemma}
For $n \in \mathbf N$ holds
\begin{eqnarray*}
\sum_{j=1}^{n}j^{2}H_{j-1}=\frac{6n(n+1)(2n+1)H_{n+1}-n(n+1)(4n+23)}{36}
\end{eqnarray*}
\end{Lemma}

{\bf Proof.}
By induction on $n$, the case $n=1$ yielding $1^{2}H_{0}=0$ on the left-hand
side and on the right-hand side
\begin{eqnarray*}
\frac{36H_{2}-54}{36}=\frac{36+18-54}{36}=0.
\end{eqnarray*}
Suppose that the equation holds for all $n\leq k$. For $n=k+1$, we have
\begin{align*}
\sum_{j=1}^{k+1}j^{2}H_{j-1}&=\sum_{j=1}^{k}j^{2}H_{j-1}+(k+1)^{2}H_{k}\\
&=\frac{6k(k+1)(2k+1)H_{k+1}-k(k+1)(4k+23)}{36}+(k+1)^{2}H_{k+1}-(k+1)\\
&=\frac{6(k+1)H_{k+1}(2k^{2}+7k+6)-k(k+1)(4k+23)-36(k+1)}{36}\\
&=\frac{6(k+1)H_{k+1}(k+2)(2k+3)-(k+1)(4k^{2}+23k+36)}{36}\\
&=\frac{6(k+1)H_{k+2}(k+2)(2k+3)-6(k+1)(2k+3)-(k+1)(4k^{2}+23k+36)}{36}\\
&=\frac{6(k+1)(k+2)H_{k+2}(2k+3)-(k+1)(4k^{2}+35k+54)}{36}\\
&=\frac{6(k+1)(k+2)H_{k+2}(2k+3)-(k+1)(k+2)(4k+27)}{36}
\end{align*}
finishing the proof of Lemma 3.3. $\bullet$
\\
\\
We also need to compute $\sum_{j=1}^{n}j(n-j+1)H_{n-j}$. We have
\begin{Lemma}
For $n \in \mathbf N$
\begin{eqnarray*}
\sum_{j=1}^{n}j(n-j+1)H_{n-j}=\frac{6nH_{n+1}(n^{2}+3n+2)
-5n^{3}-27n^{2}-22n}{36}.
\end{eqnarray*}
\end{Lemma}

{\bf Proof.} We can write $j=n+1-(n-j+1)$. Then, substituting $k=n-j+1$ we obtain
\begin{align*}
\sum_{j=1}^{n}j(n-j+1)H_{n-j}&=\sum_{j=1}^{n}[(n+1)-(n-j+1)](n-j+1)H_{n-j}\\
&=(n+1)\sum_{j=1}^{n}(n-j+1)H_{n-j}-\sum_{j=1}^{n}(n-j+1)^{2}H_{n-j}\\
&=(n+1)\sum_{k=1}^{n}kH_{k-1}-\sum_{k=1}^{n}k^{2}H_{k-1}.
\end{align*}
Thus, since we know both sums by Lemmas 3.1 and 3.3 we get
\begin{align*}
\sum_{j=1}^{n}j(n-j+1)H_{n-j}&=(n+1)\sum_{k=1}^{n}kH_{k-1}-\sum_{k=1}^{n}k^{2}H_{k-1}\\
&=(n+1)(\frac{n(n+1)H_{n+1}}{2}-\frac{n(n+5)}{4})\\
&-\frac{6n(n+1)(2n+1)H_{n+1}-n(n+1)(4n+23)}{36}\\
&=\frac{n(n+1)^{2}H_{n+1}}{2}-\frac{n(n+1)(n+5)}{4}\\
&-\frac{6n(n+1)(2n+1)H_{n+1}-n(n+1)(4n+23)}{36}\\
&=\frac{18n(n+1)^{2}H_{n+1}-9n(n+1)(n+5)}{36}\\
&-\frac{6n(n+1)(2n+1)H_{n+1}-n(n+1)(4n+23)}{36}\\
&=\frac{6nH_{n+1}(n^{2}+3n+2)-n(n+1)(5n+22)}{36}
\end{align*}
which is easily checked to be equal to the quantity in the statement above
on expanding out. $\bullet$

We are now ready to complete the evaluation of
$\sum_{j=1}^{n}M_{j-1}M_{n-j}$. Note first that
$\sum_{j=1}^{n}(n-j+1)H_{n-j}=\sum_{k=1}^{n}kH_{k-1}$ (set $k=n-j+1$) and
so Lemma 3.1 can be used to compute it. Pulling everything together, we
have
\begin{align*}
\sum_{j=1}^{n}M_{j-1}M_{n-j}&=4\sum_{j=1}^{n}jH_{j-1}(n-j+1)H_{n-j}-8n(\frac{n(n+1)H_{n+1}}{2}-\frac{n(n+5)}{4})\\
&+8(\frac{6n(n+1)(2n+1)H_{n+1}-n(n+1)(4n+23)}{36})+16n\sum_{j=1}^{n}j\\
&-8(\frac{6nH_{n+1}(n^{2}+3n+2)-5n^{3}-27n^{2}-22n}{36})-16\sum_{j=1}^{n}j^{2}\\
&+8(\frac{n(n+1)H_{n+1}}{2}-\frac{n(n+5)}{4})-16n^{2}+16\sum_{j=1}^{n}j\\
&=4\sum_{j=1}^{n}jH_{j-1}(n-j+1)H_{n-j}-4n^{2}(n+1)H_{n+1}+2n^{2}(n+5)\\
&+\frac{8}{36}[6n(n^{2}-1)H_{n+1}+n^{3}-n]+8n^{2}(n+1)-\frac{16n^{3}}{3}-8n^{2}-\frac{16n}{6}\\
&+4n(n+1)H_{n+1}-2n(n+5)-16n^{2}+8n(n+1)\\
&=4\sum_{j=1}^{n}jH_{j-1}(n-j+1)H_{n-j}-4n(n+1)(n-1)H_{n+1}\\
&+\frac{4}{3}n(n^{2}-1)H_{n+1}+\frac{176n^{3}}{36}-\frac{176n}{36}.
\end{align*}
Thus we indeed get the conclusion of Lemma 3.2, namely that
\begin{equation*}
\sum_{j=1}^{n}M_{j-1}M_{n-j}=4\sum_{j=1}^{n}jH_{j-1}(n-j+1)H_{n-j}
-\frac{8n}{3}(n^{2}-1)H_{n+1}+\frac{44n}{9}(n^{2}-1).
\end{equation*}
Returning back to the recurrence relation (1), we obtain from Lemmas 3.1 and
3.2 that
\begin{align*}
B_{n}&=\frac{(n-1)(n-2)}{2}+\frac{2(n-1)}{n}(n(n+1)H_{n+1}-\frac{5n^{2}+n}{2})\\
&+\frac{4\sum_{j=1}^{n}jH_{j-1}(n-j+1)H_{n-j}}{n}-\frac{8}{3}(n^{2}-1)H_{n+1}+\frac{44}{9}(n^{2}-1)+\frac{2}{n}\sum_{j=1}^{n}B_{j-1}\\
&=\frac{(n-1)(n-2)}{2}+2(n-1)(n+1)H_{n+1}-(n-1)(5n+1)\\
&+\frac{4\sum_{j=1}^{n}jH_{j-1}(n-j+1)H_{n-j}}{n}-\frac{8}{3}(n^{2}-1)H_{n+1}
+\frac{44}{9}(n^{2}-1)+\frac{2}{n}\sum_{j=1}^{n}B_{j-1}.
\end{align*}
Finally,
\begin{align*}
B_{n}&=\frac{4\sum_{j=1}^{n}jH_{j-1}(n-j+1)H_{n-j}}{n}+\frac{2}{n}\sum_{j=1}^{n}B_{j-1}+\frac{-9n^{2}+5n+4}{2}\\
&-\frac{2}{3}(n^{2}-1)H_{n+1}+\frac{44}{9}(n^{2}-1).
\end{align*}
Multiplying by $n$, we have
\begin{align*}
nB_{n}&=4\sum_{j=1}^{n}jH_{j-1}(n-j+1)H_{n-j}+2\sum_{j=1}^{n}B_{j-1}+n\frac{-9n^{2}+5n+4}{2}\\
&-n\frac{2}{3}(n^{2}-1)H_{n+1}+n\frac{44}{9}(n^{2}-1).
\end{align*}
For $n+1$, we have similarly
\begin{align*}
&(n+1)B_{n+1}\\
&=4\sum_{j=1}^{n+1}jH_{j-1}(n-j+2)H_{n+1-j}
+2\sum_{j=1}^{n+1}B_{j-1}
+(n+1)\frac{-9(n+1)^{2}+5(n+1)+4}{2}\\
&-(n+1)\frac{2}{3}[(n+1)^{2}-1]H_{n+2}+(n+1)\frac{44}{9}[(n+1)^{2}-1].
\end{align*}
Subtracting $nB_{n}$ from $(n+1)B_{n+1}$, we obtain
\begin{align*}
&(n+1)B_{n+1}-nB_{n}\\
&=4[\sum_{j=1}^{n+1}jH_{j-1}(n-j+2)H_{n+1-j}-\sum_{j=1}^{n}jH_{j-1}(n-j+1)H_{n-j}]\\
&+2B_{n}+(n+1)\frac{-9(n+1)^{2}+5(n+1)+4}{2}-n\frac{-9n^{2}+5n+4}{2}\\
&-(n+1)\frac{2}{3}[(n+1)^{2}-1]H_{n+2}+\frac{2}{3}n(n^{2}-1)H_{n+1}\\
&+(n+1)\frac{44}{9}[(n+1)^{2}-1]-n\frac{44}{9}(n^{2}-1)\\
&=4[\sum_{j=1}^{n}jH_{j-1}(n-j+2)H_{n+1-j}-\sum_{j=1}^{n}jH_{j-1}(n-j+1)H_{n-j}]\\
&+2B_{n}-\frac{27n^{2}+17n}{2}-\frac{2}{3}n(n+1)(n+2)H_{n+2}
+\frac{2}{3}nH_{n+1}(n^{2}-1)+\frac{44n(n+1)}{3}
\end{align*}
noting that the term for $j=n+1$ does not contribute to the sum.
In the first sum, we use the facts that
$H_{n+1-j}=H_{n-j}+1/(n+1-j)$ and that $n-j+2=(n-j+1)+1$, and then we get
\begin{align*}
&4[\sum_{j=1}^{n}jH_{j-1}(n-j+1)H_{n-j}+\sum_{j=1}^{n}\frac{jH_{j-1}(n-j+1)}{n-j+1}\\
&+\sum_{j=1}^{n}jH_{j-1}H_{n-j+1}-\sum_{j=1}^{n}jH_{j-1}(n-j+1)H_{n-j}]+2B_{n}\\
&-\frac{27n^{2}+17n}{2}-\frac{2}{3}n(n+1)(n+2)H_{n+2}
+\frac{2}{3}nH_{n+1}(n^{2}-1)+\frac{44n(n+1)}{3}.
\end{align*}
The first sum on the first line cancels with the equal
sum on the second line, the second sum on the first line simplifies, and again using $H_{n+2}=H_{n+1}+1/(n+2)$ on the last line, we obtain
\begin{align*}
&4[\sum_{j=1}^{n}jH_{j-1}+\sum_{j=1}^{n}jH_{j-1}H_{n-j+1}]
+2B_{n}-2n(n+1)H_{n+1}+\frac{1}{2}n(n+11).
\end{align*}
We thus see that we have to work out the following expression:
\begin{eqnarray*}
\sum_{j=1}^{n}jH_{j-1}H_{n+1-j}.
\end{eqnarray*}
We note that

\begin{align*}
\sum_{j=1}^{n}jH_{j-1}H_{n+1-j}&=1H_{0}H_{n}+2H_{1}H_{n-1}+3H_{2}H_{n-2}+\ldots+(n-1)H_{n-2}H_{2}+nH_{n-1}H_{1}\\
&=\frac{n+2}{2}(H_{1}H_{n-1}+H_{2}H_{n-2}+\ldots+H_{n-2}H_{2}+H_{n-1}H_{1}+H_{n}H_{0})\\
&=\frac{n+2}{2}\sum_{j=1}^{n}H_{j}H_{n-j}
\end{align*}
so it suffices now to obtain the quantity $\sum_{j=1}^{n}H_{j}H_{n-j}$.
\\
\\
Sedgewick \cite{sedgewick} presents and proves the following result:
\begin{Lemma}
\begin{align*}
\sum_{i=1}^{n}H_{i}H_{n+1-i}=
   (n+2)(H^{2}_{n+1}-H^{(2)}_{n+1})-2(n+1)(H_{n+1}-1).
\end{align*}
\end{Lemma}

{\bf Proof.}
\begin{align*}
\sum_{i=1}^{n}H_{i}H_{n+1-i}
=\sum_{i=1}^{n}H_{i}\sum_{j=1}^{n+1-i}\frac{1}{j}.
\end{align*}
The set $\{(i,j):\,1\leq i\leq n, \,1\leq j\leq n+1-i\}$ is, as a picture
easily shows, the same as $\{(i,j):\,1\leq i\leq n+1-j,\,1\leq j\leq n\}$.
Thus the above is
\begin{align*}
&\sum_{j=1}^{n}\frac{1}{j}\sum_{i=1}^{n+1-j}H_{i}
=\sum_{j=1}^{n}\frac{1}{j}[(n+2-j)H_{n+1-j}-(n+1-j)]
\end{align*}
To see the claim about the sum of the $H_{j}$s, we note that
\begin{align*}
&\sum_{j=1}^{n}H_{j}=H_{1}+H_{2}+\ldots+H_{n}=1+(1+\frac{1}{2})+\ldots+(1+\frac{1}{2}+\ldots+\frac{1}{n})\\
&=n+(n-1)\frac{1}{2}+\ldots+[n-(n-1)]\frac{1}{n}=n(1+\frac{1}{2}+\ldots+\frac{1}{n})-(\frac{1}{2}+\frac{2}{3}+\ldots+\frac{n-1}{n})\\
&=n(1+\frac{1}{2}+\ldots+\frac{1}{n})-[(1-\frac{1}{2})+(1-\frac{1}{3})+\ldots+(1-\frac{1}{n})]\\
&=nH_{n}-n+H_{n}=(n+1)H_{n}-n
\end{align*}
Thus we get, reusing the result about $\sum_{j=1}^{n}H_{j}$,  
\begin{align*}
&(n+2)\sum_{j=1}^{n}\frac{H_{n+1-j}}{j}-\sum_{j=1}^{n}H_{n+1-j}-(n+1)H_{n}+n\\
&=(n+2)\sum_{j=1}^{n}\frac{H_{n+1-j}}{j}-[(n+1)H_{n}-n]-(n+1)H_{n}+n\\
&=(n+2)\sum_{j=1}^{n}\frac{H_{n+1-j}}{j}-2(n+1)(H_{n+1}-1).
\end{align*}
To analyse the first sum above, we note (again following \cite{sedgewick} here)
\begin{align*}
\sum_{j=1}^{n}\frac{H_{n+1-j}}{j}&=\sum_{j=1}^{n}\frac{H_{n-j}}{j}+\sum_{j=1}^{n}\frac{1}{j(n+1-j)}\\
&=\sum_{j=1}^{n-1}\frac{H_{n-j}}{j}+\frac{1}{n+1}\sum_{j=1}^{n}(\frac{1}{j}+\frac{1}{n+1-j})
\end{align*}
and this gives that
\begin{align}
&\sum_{j=1}^{n}\frac{H_{n+1-j}}{j}=\sum_{j=1}^{n-1}
\frac{H_{n-j}}{j}+2\frac{H_{n}}{n+1}.
\end{align}
Iterating this equation, and using $H_{0}=0$ at the end, we obtain the identity
\begin{eqnarray}
\sum_{j=1}^{n}\frac{H_{n+1-j}}{j}=2\sum_{k=1}^{n}\frac{H_{k}}{k+1}.
\end{eqnarray}
The right-hand side is
\begin{align*}
&2\sum_{k=1}^{n}\frac{H_{k}}{k+1}=2\sum_{k=2}^{n+1}\frac{H_{k-1}}{k}
=2\sum_{k=1}^{n+1}\frac{H_{k-1}}{k}=2\sum_{k=1}^{n+1}\frac{H_{k}}{k}
-2\sum_{k=1}^{n+1}\frac{1}{k^{2}}\\
&=2\sum_{k=1}^{n+1}\sum_{j=1}^{k}\frac{1}{jk}-2H^{(2)}_{n+1}
=2\sum_{j=1}^{n+1}\sum_{k=j}^{n+1}\frac{1}{jk}-2H^{(2)}_{n+1}
=2\sum_{k=1}^{n+1}\sum_{j=k}^{n+1}\frac{1}{kj}-2H^{(2)}_{n+1}.
\end{align*}
Again noting that $\{(j,k):\,k\leq j\leq n+1,\,1\leq k\leq n+1\}$ gives
the same values of $1/(jk)$ as $\{(j,k): 1\leq j\leq k,\,1\leq k\leq n+1\}$
provided we note that the terms for $j=k$ are repeated, we get
\begin{align*}
&2\sum_{k=1}^{n+1}\sum_{j=k}^{n+1}\frac{1}{kj}-2H^{(2)}_{n+1}
=\sum_{k=1}^{n+1}(\sum_{j=1}^{n+1}\frac{1}{kj}+\frac{1}{k^{2}})-2H^{(2)}_{n+1}\\
&=H^{2}_{n+1}+H^{(2)}_{n+1}-2H^{(2)}_{n+1}
=H^{2}_{n+1}-H^{(2)}_{n+1}.
\end{align*}
Thus, we have, as in the statement of Lemma 3.5, 
\begin{eqnarray*}
\sum_{i=1}^{n}H_{i}H_{n+1-i}
=(n+2)(H^{2}_{n+1}-H^{(2)}_{n+1})-2(n+1)(H_{n+1}-1).
\end{eqnarray*}

Also, the following Corollary is obtained, using equations from the
last Lemma, by Sedgewick \cite{sedgewick}.
\begin{Corollary} For $n \in \mathbf N$, it holds
\begin{align*}
H_{n+1}^{2}-H_{n+1}^{(2)}=2\sum_{j=1}^{n}\frac{H_{j}}{j+1}
\end{align*}
\end{Corollary}
{\bf Proof.} From equations (3) and (4), we see that
\begin{align*}
&H^{2}_{n+1}-H^{(2)}_{n+1}=\sum_{j=1}^{n}\frac{H_{n+1-j}}{j}
=\sum_{j=1}^{n-1}\frac{H_{n-j}}{j}+2\frac{H_{n}}{n+1}\\
&\Rightarrow H^{2}_{n+1}-H^{(2)}_{n+1}=H^{2}_{n}-H^{(2)}_{n}
+2\frac{H_{n}}{n+1}\\
&\Rightarrow H^{2}_{n+1}-H^{(2)}_{n+1}=2\sum_{j=1}^{n}\frac{H_{j}}{j+1}
\end{align*}
by iteration. $\bullet$

We will use the above Lemma and Corollary in our analysis. We have that
\begin{eqnarray*}
\sum_{i=1}^{n}H_{i}H_{n+1-i}=\sum_{i=1}^{n}[H_{i}(H_{n-i}+\frac{1}{n+1-i})]=\sum_{i=1}^{n}H_{i}H_{n-i}+\sum_{i=1}^{n}\frac{H_{i}}{n+1-i}
\end{eqnarray*}
The second sum substituting $j=n+1-i$ becomes
\begin{eqnarray*}
\sum_{i=1}^{n}\frac{H_{i}}{n+1-i}=\sum_{j=1}^{n}\frac{H_{n+1-j}}{j}
\end{eqnarray*}
As we have seen this is equal to
\begin{eqnarray*}
\sum_{j=1}^{n}\frac{H_{n+1-j}}{j}=H^{2}_{n+1}-H^{(2)}_{n+1}
\end{eqnarray*}
Hence, by Lemma 3.5
\begin{align*}
\sum_{i=1}^{n}H_{i}H_{n-i}&=(n+2)(H^{2}_{n+1}-H^{(2)}_{n+1})-2(n+1)(H_{n+1}-1)-(H^{2}_{n+1}-H^{(2)}_{n+1})\\
&=(n+1)[(H^{2}_{n+1}-H^{(2)}_{n+1})-2(H_{n+1}-1)].
\end{align*}
Using the above equation and the result obtained in page 13, just before Lemma
3.5, we deduce that
\begin{align*}
\sum_{j=1}^{n}jH_{j-1}H_{n+1-j}&=\frac{n+2}{2}(n+1)[(H^{2}_{n+1}-H^{(2)}_{n+1})-2(H_{n+1}-1)]\\
&={n+2\choose 2}[(H^{2}_{n+1}-H^{(2)}_{n+1}-2(H_{n+1}-1)].
\end{align*}

Having worked out all the expressions involved in the following relation,
we can now finish off: 
\begin{align*}
(n+1)B_{n+1}-nB_{n}=
4(\sum_{j=1}^{n}jH_{j-1}+\sum_{j=1}^{n}jH_{j-1}H_{n-j+1})+2B_{n}-2n(n+1)H_{n+1}+\frac{1}{2}n(n+11).
\end{align*}
We have
\begin{align*}
&(n+1)B_{n+1}-nB_{n}\\
&=4[\frac{n(n+1)H_{n+1}}{2}-\frac{n(n+5)}{4}+{n+2\choose 2}[(H^{2}_{n+1}-H^{(2)}_{n+1})-2(H_{n+1}-1)]]\\
&+2B_{n}-2n(n+1)H_{n+1}+\frac{1}{2}n(n+11)\\
&=2n(n+1)H_{n+1}-n(n+5)+2(n+1)(n+2)(H^{2}_{n+1}-H^{(2)}_{n+1})\\
&-4(n+1)(n+2)(H_{n+1}-1)+2B_{n}-2n(n+1)H_{n+1}+\frac{1}{2}n(n+11)\\
&=2(n+1)(n+2)(H^{2}_{n+1}-H^{(2)}_{n+1})-4(n+1)(n+2)(H_{n+1}-1)
-\frac{n(n-1)}{2}+2B_{n}.
\end{align*}
Then,
\begin{align*}
B_{n+1}&=2(n+2)(H^{2}_{n+1}-H^{(2)}_{n+1})-4(n+2)(H_{n+1}-1)-\frac{n(n-1)}{2(n+1)}+\frac{n+2}{n+1}B_{n}\\
\frac{B_{n+1}}{n+2}&=\frac{B_{n}}{n+1}+2(H^{2}_{n+1}-H^{(2)}_{n+1})
-4(H_{n+1}-1)-\frac{n(n-1)}{2(n+1)(n+2)}.
\end{align*}
The last equation is equivalent to
\begin{align*}
\frac{B_{n}}{n+1}=\frac{B_{n-1}}{n}+2(H^{2}_{n}-H^{(2)}_{n})-4(H_{n}-1)
-\frac{(n-1)(n-2)}{2n(n+1)}.
\end{align*}
Iterating the recurrence relation, we obtain
\begin{eqnarray*}
\frac{B_{n}}{n+1}=\frac{B_{0}}{1}+2\sum_{i=1}^{n}(H^{2}_{i}-H^{(2)}_{i})
-4\sum_{i=1}^{n}(H_{i}-1)-\sum_{i=1}^{n}\frac{(i-1)(i-2)}{2i(i+1)}.
\end{eqnarray*}
Since $B_{0}=0$, it is
\begin{eqnarray*}
\frac{B_{n}}{n+1}=2\sum_{i=1}^{n}(H^{2}_{i}-H^{(2)}_{i})
-4\sum_{i=1}^{n}(H_{i}-1)-\sum_{i=1}^{n}\frac{(i-1)(i-2)}{2i(i+1)}.
\end{eqnarray*}
The first sum, by Corollary 3.6 is equal to
\begin{align*}
\sum_{i=1}^{n}(H^{2}_{i}-H^{(2)}_{i})
&=(H_{1}^{2}-H_{1}^{(2)})+(H_{2}^{2}-H_{2}^{(2)})+\ldots +(H_{n}^{2}-H_{n}^{(2)})\\
&=n(H^{2}_{n}-H^{(2)}_{n})-2\sum_{i=1}^{n-1}\frac{iH_{i}}{i+1}\\
&=(n+1)(H^{2}_{n}-H^{(2)}_{n})-2\sum_{i=1}^{n-1}\frac{iH_{i}}{i+1}-2\sum_{i=1}^{n-1}\frac{H_{i}}{i+1}\\
&=(n+1)(H^{2}_{n}-H^{(2)}_{n})-2\sum_{i=1}^{n-1}H_{i}\\
&=(n+1)(H^{2}_{n}-H^{(2)}_{n})+2n-2nH_{n}.
\end{align*}
Note that on the third line we add and subtract simultaneously $(H^{2}_{n}-H^{(2)}_{n})$, which is equal to
$2\sum_{i=1}^{n-1}H_{i}/i+1$ by Corollary 3.6. Doing so, the fraction is cancelled and
the corresponding sum can be easily computed. Hence
\begin{align*}
\frac{B_{n}}{n+1}&=2(n+1)(H^{2}_{n}-H^{(2)}_{n})+4n-4nH_{n}-4[(n+1)H_{n}-2n]-\sum_{i=1}^{n}\frac{(i-1)(i-2)}{2i(i+1)}\\
&=2(n+1)(H^{2}_{n}-H^{(2)}_{n})+4n-4nH_{n}-4[(n+1)H_{n}-2n]-\sum_{i=1}^{n}(\frac{i+2}{2i}-\frac{3}{i+1})\\
&=2(n+1)(H^{2}_{n}-H^{(2)}_{n})+12n-8nH_{n}-4H_{n}-\sum_{i=1}^{n}\frac{i+2}{2i}+\sum_{i=1}^{n}\frac{3}{i+1}\\
&=2(n+1)(H^{2}_{n}-H^{(2)}_{n})+12n-8nH_{n}-4H_{n}-\frac{n}{2}-H_{n}+3H_{n+1}-3\\
&=2(n+1)(H^{2}_{n}-H^{(2)}_{n})+12n-8nH_{n}-4H_{n}-\frac{n}{2}-H_{n}+3H_{n}+\frac{3}{n+1}-3\\
&=2(n+1)(H^{2}_{n}-H^{(2)}_{n})-H_{n}(8n+2)+\frac{23n}{2}+\frac{3}{n+1}-3\\
&=2(n+1)(H^{2}_{n}-H^{(2)}_{n})-H_{n}(8n+2)+\frac{23n^{2}+17n}{2(n+1)}.
\end{align*}
Finally, multiplying both sides by $n+1$ we obtain
\begin{eqnarray*}
B_{n}=2(n+1)^{2}(H^{2}_{n}-H^{(2)}_{n})
-H_{n}(n+1)(8n+2)+\frac{23n^{2}+17n}{2}.
\end{eqnarray*}
Now, the Proof of Theorem 2.3 is complete. Consequently, the Variance of the number of pairwise comparisons $C_{n}$ of Randomised Quicksort is equal to
\begin{eqnarray*}
7n^{2}-4(n+1)^{2}H_{n}^{(2)}-2(n+1)H_{n}+13n.
\end{eqnarray*}
\pagebreak

\end{document}